\documentclass{amsart}
\usepackage{amsthm,amsmath,amsfonts}
\usepackage{mathrsfs}
\usepackage[colorlinks=false,linktocpage=true]{hyperref}
\usepackage{tocvsec2}
\usepackage{latexsym}
\numberwithin{equation}{section}
\font\tenscrpt=eusm10
\font\sevenscrpt=eusm10 scaled 700
\font\fivescrpt=eusm10 scaled 500
\newfam\eusmfam
\textfont\eusmfam=\tenscrpt
\scriptfont\eusmfam=\sevenscrpt
\scriptscriptfont\eusmfam=\fivescrpt
\def\scr#1{{\fam\eusmfam\relax#1}}

\numberwithin{equation}{section}
\font\tenscrpt=eusm10
\font\sevenscrpt=eusm10 scaled 700
\font\fivescrpt=eusm10 scaled 500
\newfam\eusmfam
\textfont\eusmfam=\tenscrpt
\scriptfont\eusmfam=\sevenscrpt
\scriptscriptfont\eusmfam=\fivescrpt
\def\scr#1{{\fam\eusmfam\relax#1}}
\def\qed{\quad\vcenter{\hrule\hbox{\vrule height.6em\kern.6em\vrule}\hrule}}

\newenvironment{pf*}[1]{{\textsc #1.\quad}}{$\qed$\bigskip\newline}
\newtheorem{thm}{Theorem}[section]
\newtheorem{cor}{Corollary}[section]

\newtheorem{rem}{Remark}[section]

\newcommand{\thmref}[1]{Theorem~{\rm $\ref{#1}$}}

\newcommand{\corref}[1]{Corollary~{\rm $\ref{#1}$}}

\newcommand{\remref}[1]{Remark~{\rm $\ref{#1}$}}
\newcommand{\eqnref}[1]{{\rm (\ref{#1})}}
\def\eqdef{\overset{\triangle}{=}}
\def\half{\frac12}

\def\taunU{\tau_n^U}
\def\taunV{\tau_n^V}
\def\ehNabh{e_{heat}^{Neu}(a,b,h)}
\def\ehNabdh{e_{heat}^{Neu}(a,b+d,h)}
\def\ehNazh{e_{heat}^{Neu}(a,0,h)}
\def\ehabh{e_{heat}(a,b,h)}
\def\ehabdh{e_{heat}(a,b+d,h)}

\def\P{{\mathbb P}}
\def\Q{{\mathbb Q}}
\def\Ptn{\tilde{\P}_n}
\def\L{\mathbb L}
\def\LawXP{\L_{\P}^{X}}
\def\LawUP{\L_{\P}^{U}}
\def\LawVQ{\L_{\Q}^{V}}
\def\LawUtP{\L_{\P}^{U(t,\cdot)}}
\def\LawVtQ{\L_{\Q}^{V(t,\cdot)}}
\def\LawUtxP{\L_{\P}^{U(t,x)}}
\def\LawVtxQ{\L_{\Q}^{V(t,x)}}
\def\EP{{\mathbb E}_{\P}}

\def\Ru{R_u}
\def\RU{R_U}
\def\RUsq{\RU^2}
\def\RV{R_V}
\def\RNfTtaunW{\Xi_{T\wedge\taunU}^{\RU,{\Wm}}(\R)}

\def\Wm{\scr W}
\def\Wtm{\tilde{\Wm}}
\def\WtB{\Wm_t^{(2)}(B)}
\def\WttB{\Wtm_t(B)}
\def\Wtmn{\Wtm_n}
\def\WttBn{\Wtm_{t\wedge\taunU}(B)}
\def\solUWt{(U,\Wtmn)}

\def\Oo{\Omega^{(1)}}
\def\Ot{\Omega^{(2)}}

\def\BR{{\scr B}(\R)}
\def\N{{\mathbb N}}
\def\R{{\mathbb R}}

\def\RTL{{\scr R}_{T,L}}

\def\CRTLR{C(\RTL;\R)}

\def\RtaunT{{\scr R}_{T\wedge\taunU}}

\def\RT{{\scr R}_{T}}

\def\RTTcl{{\overset{\smile}{\scr R}}_{T}}
\def\CRTR{C(\RT;\R)}

\def\filpspace{(\Ot, \scr{F}, \{{\scr{F}}_t\},\P)}
\def\filqspace{(\Oo, \scr{H}, \{{\scr{H}}_t\},\Q)}
\def\filTpspace{(\Ot, \scr{F}_T, \{{\scr{F}}_t\},\Ptn)}

\def\Ft{\scr{F}_t}
\def\FT{\scr{F}_T}
\begin{document}
\title[Law equivalence and compact support for SPDEs]{SPDEs law  equivalence and the compact support
property: applications to the Allen-Cahn SPDE}
\author{Hassan ALLOUBA} 
\address{Department of Mathematics and Statistics, University of Massachusetts, Amherst, MA 01003-4515}
\email{allouba@math.umass.edu}
\date{May 1, 2000}
\begin{abstract}
Using our uniqueness in law transfer result for SPDEs, described in a recent note, we prove the equivalence of
laws of SPDEs differing by a drift, under vastly applicable conditions.  This gives us the equivalence
in the compact support property among a large class of SPDEs.  As an
important application, we prove the equivalence in law of the Allen-Cahn and the associated heat SPDEs; and
we give a criterion for the compact support property to hold for the
Allen-Cahn SPDE with diffusion function $a(t,x,u)=Cu^\gamma$, with $C\ne0$ and
$ 1/2\le\gamma<1$.
\end{abstract}
\maketitle
\section{Statements and discussions of results.}
We start by considering the pair of parabolic SPDEs 
\begin{equation}
 \begin{cases} \displaystyle\frac{\partial U}{\partial t}=\Delta _{x}U+b(t,x,U)+a(t,x,U)
\displaystyle\frac{\partial^2 W}{\partial t\partial x}; & (t,x)\in\RTTcl, \cr U_x(t,-\infty)=U_x(t,\infty)=0; & 0<t\le T,
\cr U(0,x)=h(x); & x\in\R, \end{cases} \label{Pb} \end{equation} 
and
\begin{equation} \begin{cases}
\displaystyle\frac{\partial V}{\partial t}=\Delta _{x}V+(b+d)(t,x,V)+a(t,x,V) \displaystyle\frac{\partial^2 W}{\partial
t\partial x}; & (t,x)\in\RTTcl, \cr V_x(t,-\infty)=V_x(t,\infty)=0; & 0<t\le T, \cr V(0,x)=h(x); & x\in\R.
\end{cases}
\label{Pbd} \end{equation}
on $\RT\eqdef[0,T]\times\R$, where  $W(t,x)$ is
the Brownian sheet corresponding to the driving space-time white noise,
written formally as $\partial^2W/\partial t\partial x$.  As in Walsh \cite{WA},
white noise is regarded as a continuous orthogonal martingale
measure, which we denote by $\Wm$.  The diffusion
$a(t,x,u)$ and the drifts $b(t,x,u)$ and $d(t,x,u)$ are Borel-measurable $\R$-valued functions on $\RT\times\R$; and
$h:\R\to\R$ is a bounded continuous function.
Henceforth, we will denote \eqnref{Pb} and
\eqnref{Pbd} by $\ehNabh$ and $\ehNabdh$, respectively.  When $b\equiv0$, we denote
\eqnref{Pb} by $\ehNazh$.  In the interest of getting quickly to our main results,
we refer the reader to \cite{Allouba1} for the rigorous interpretation of all SPDEs considered in this paper, with the obvious
modifications to accomodate the change of space from $\RTL=[0,T]\times[0,L]$ to $\RT$.  Also, the law of a random variable $X$
under the probability measure $\P$ is denoted by $\LawXP$.   Proceeding toward a precise statement of our results, let
$\Ru(t,x)\eqdef d(t,x,u)/a(t,x,u)$, for any $(t,x,u)\in\RT\times\R$, whenever the ratio is well defined.
Let $\lambda$ denote Lebesgue measure.  Our law equivalence result for the pair $\ehNabh$ and
$\ehNabdh$ can now be stated as
\begin{thm}
Let $(V,W^{(1)})$ be a solution $($weak or strong\/$)$ to $\ehNabdh$ on some probability space
$\filqspace$.  Assume
that $\RU$ and $\RV$ are in $L^2(\RT,\lambda)$,
almost surely, whenever the random fields $U$ and $V$ solve $($weakly or strongly\/$)$ $\ehNabh$ and $\ehNabdh$,
respectively.  Assume further that there is a unique-in-law solution $(U,W^{(2)})$ to the heat SPDE $\ehabh$ on $\filpspace$.
Then $\LawVQ$  and $\LawUP$ are mutually absolutely continuous $($\!on $(\CRTR)$\/$)$.
\label{AC1}
\end{thm}
\begin{rem}
By Theorem 1.1 in \cite{Allouba1}, which can trivially be extended from $\RTL$ to $\RT$ $($replacing $[0,L]$ with $\R$ in the $L^2$ condition
and in the proof$),$ uniqueness
in law for the SPDE $\ehabh$ is equivalent to uniqueness in law for $\ehabdh$  under
the $L^2(\RT,\lambda)$ condition on $\RU$ and $\RV$.  So, we can replace the uniqueness assumption on $\ehabh$ in
\thmref{AC1} above by that on $\ehabdh$.   Also, our Neumann conditions may be replaced by Dirichlet conditions without
affecting the conclusion of \thmref{AC1}.  Finally, \thmref{AC1} and its proof are valid when $\RT$ is
replaced by $\RTL$ $($replacing $\R$ with $[0,L]$$)$.
\label{equiv}
\end{rem}
As an immediate consequence of \thmref{AC1} and \remref{equiv},
we get the following law equivalence between the Allen-Cahn SPDE 
\begin{equation}
\begin{cases}
\displaystyle\frac{\partial V}{\partial t}=\Delta _{x}V+2V(1-{V}^{2})+ C V^\gamma
\displaystyle\frac{\partial^2 W}{\partial t\partial x};
& (t,x)\in\RTTcl,\ C>0,
\cr V_x(t,-\infty)=V_x(t,\infty)=0; & 0<t\le T,
\cr V(0,x)=h(x); & x\in\R,
\end{cases}
\label{AC}
\end{equation}and its associated heat SPDE
(the one obtained from the Allen-Cahn SPDE \eqref{AC} by removing the Allen-Cahn nonlinearity).  We note here that the proof of the
uniqueness for the Allen-Cahn SPDE in Theorem 1.2 in \cite{Allouba1} works just as well for the case $\gamma=1/2$, in addition
to $\frac12<\gamma<1$, because
the SPDE in \eqnref{Pb} with $b\equiv0$ and $a(t,x,u)=Cu^{1/2}$ admits uniqueness in law as discussed in \cite{Mu} p.~326 and in
\cite{RC}.
\begin{cor}
Suppose that $V$ and $U$ solve $($weakly or strongly\/$)$ the Allen-Cahn SPDE and its
associated heat SPDE, respectively, on $\RTL$ $($see $(0.3)$ in \cite{Allouba1})  and with $1/2\le\gamma\le1$.   Then the
laws of $U$ and $V$  are equivalent $($\!on $(\CRTLR)$\/$)$.  If $\RTL$ is replaced with $\RT,$ if $1/2\le\gamma<1,$ if
$h(x)$ has compact support, and if $\RV$ is in $L^2(\RT,\lambda)$ a.s.$:$
\begin{equation}
\int_{\RT}\RV^2(t,x)dtdx=\frac{4}{C^2}\int_{\RT}V^{2(1-\gamma)}(V^4-2V^2+1)dtdx<\infty; \ \mbox{almost surely},
\label{integ}
\end{equation}
then the laws of $U$ and $V$  are equivalent $($\!on $(\CRTR)$\/$)$.
\label{ACandheat}
\end{cor}
\begin{rem}
In the first part of \corref{ACandheat}, the continuity of $U$ and $V$ insures that $\RU$ and $\RV$ are in $L^2(\RTL,\lambda)$,
for $0\le\gamma\le1$ $($see the proof of Theorem 1.2 in \cite{Allouba1}\/$)$.
When $\RTL$ is replaced by $\RT$, we do not require that $\RU$ be in
$L^2(\RT,\lambda)$ $(\eqref{integ}$ with $U$ instead of $V$$)$.  This is because $\RU$ is already in $L^2(\RT,\lambda),$ since
$U(t,\cdot)$ has compact support for each $t$
in the range $1/2\le\gamma<1$ by \cite{Mu,Kr}.  Also, when we replace $\RTL$ with $\RT$, $\gamma\le1$ is replaced with
$\gamma<1;$ since, in this case, when $\gamma=1$ the integrability assumption in $\eqref{integ}$ has obvious problems for both $U$ and $V$.
\end{rem}
 Let $\Omega=\CRTR$ and denote elements of $\Omega$ by $\omega$.  Let $X$ be the coordinate mapping process on
 $\Omega$: $X_\omega(t,x)\eqdef\omega(t,x)$.  Denote by ${\scr G_{t,x}^X}$,
 ${\scr G_{t,\cdot}^X}$, and $\scr G_{\cdot,\cdot}^X$
the sigma fields of subsets of $\Omega$ generated by $X$  when $(t,x)$ is fixed, when $t$ is fixed but $x$ is not,
and when both $t$ and $x$ are not fixed, respectively.  I.e.,
\begin{align*}
&{\scr G_{t,x}^X}=\sigma\Big(\big\{\omega\in\Omega;X_\omega(t,x)=\omega(t,x)\in A\big\}; \ A\in\scr B(\R)\Big);\quad(t,x)\in\RT,\\
&{\scr G_{t,\cdot}^X}=\sigma\Big(\big\{\omega\in\Omega;\big(X_\omega(t,x_1)=\omega(t,x_1),\ldots,
X_\omega(t,x_n)=\omega(t,x_n)\big)\in A\big\};\Big. &\\ &\Big. \hspace{1.75in}n\ge1,\, A\in\scr B(\R^n),\, x_i\in\R,\, i=1,\ldots,n\Big);
\quad t\in[0,T],\\
&{\scr G_{\cdot,\cdot}^X}=\sigma\Big(\big\{\omega\in\Omega;\big(X_\omega(t_1,x_1)=\omega(t_1,x_1),\ldots,
X_\omega(t_n,x_n)=\omega(t_n,x_n)\big)\in A\big\};\Big. &\\ &\Big.\hspace{2.15in}  n\ge1,\, A\in\scr B(\R^n),\, (t_i,x_i)\in\RT,\, i=1,\ldots,n\Big).
\end{align*}
Then, clearly, ${\scr G_{t,x}^X}\subseteq{\scr G_{t,\cdot}^X}\subseteq{\scr B(\CRTR)}=\scr G_{\cdot,\cdot}^X$,  the last
equality is a trivial extension of Problem 4.2 p.~60 in \cite{Shreve},
and so absolute continuity on ${\scr B(\CRTR)}$ implies absolute continuity on ${\scr G_{t,x}^X}$ and
${\scr G_{t,\cdot}^X}$.
This observation along with \thmref{AC1} easily give us
\begin{cor} Under the conditions of \thmref{AC1}, $\LawVtxQ$ is equivalent to $\LawUtxP$ on $\R,$ for every $(t,x)\in\RT$
$($in particular, if one is absolutely continuous with respect to Lebesgue measure then so is the other\/$);$ and $\LawVtQ$
is equivalent to $\LawUtP$ on $C(\R;\R),$ for every $t\in[0,T]$.
\label{AC2}
\end{cor}
By proving law equivalence between $\ehNabh$ and $\ehNabdh$ under considerably weaker conditions,
\thmref{AC1} and \corref{AC2} extend and make more applicable the notion of relative absolute
continuity in our earlier work (Theorem 3.3.3 in \cite{HAd} or Theorem 4.3 in \cite{HA1}).   Like Theorem 4.3 in \cite{HA1},
\thmref{AC1} and \corref{AC2} (and thus the first assertion of \thmref{CSPU} below) are equally valid for wave SPDEs, space-time SDEs, and SDEs (cf. Theorem 3.7, Theorem 4.3,
Theorem 5.3, and their proofs under the stronger conditions of \cite{HA1}).  An interesting application of
\thmref{AC1} and \corref{AC2} is to allow us to prove the following
theorem about the compact support property of solutions to a large class of SPDEs containing the Allen-Cahn SPDE:
\begin{thm}
Assume that the conditions of \thmref{AC1} hold.  Then, $U(\cdot,\cdot)$ $(U(t,\cdot))$ has compact support
iff $V(\cdot,\cdot)$ $(V(t,\cdot))$ does.  In particular, if $V(t,x)$ is a solution to the Allen-Cahn SPDE \eqref{AC}, $h(x)$ has compact
support, and $\half\le\gamma<1;$ then, for each $t\in[0,T],$ $V(t,\cdot)$ has compact support as a function of $x$ iff \eqref{integ}
holds.
\label{CSPU}
\end{thm}
It is noteworthy that all the Allen-Cahn SPDE results and their proofs here and in \cite{Allouba1} are valid for the
KPP SPDE, obtained by replacing the Allen-Cahn term $2V(1-V^2)$ by the KPP term $V(1-V)$.   In \cite{HAd,HA2}, we gave
a proof of the existence of solutions to heat SPDEs with continuous diffusion coefficient $a$ and measurable
drift $b$---with $a$ satisfying a linear growth condition and $b/a$ satisfying Novikov's condition---using a system of stochastic
differential-difference equations (SDDEs).   In \cite{HA},  we use our SDDE approach and the results of this note and \cite{Allouba1} to
further investigate the existence and some properties of SPDEs considered here.
 \section{Proofs of results}
\begin{pf*}{Proof of \thmref{AC1}}
It follows  from the uniqueness in law assumption for $\ehNabh$,
the almost sure $L^2(\RT,\lambda)$ condition on $\RV$, and a trivial extension of Theorem 1.1 in \cite{Allouba1} to the space $\RT$
that we have  uniqueness in law for $\ehNabdh$ (see \remref{equiv}).

Now, take $\{\taunU\}$ and $\{\taunV\}$ to be the sequences of stopping times
\begin{equation}
\taunU\eqdef T\wedge\inf\left\{0\le t\le T;\intop_{[0,t]\times\R}\RUsq(s,x)ds dx=n\right\};
\ n\in\N,
\label{STimes}
\end{equation}
and $\{\taunV\}$ is gotten from \eqref{STimes} by replacing $U$ with $V$.  Let $\Wtm=\{\WttB,\Ft;0\le t\le T,B\in\BR\}$
be given by
$$\WttB\eqdef\WtB-\int _{[0,t]\times B}\RU(s,x)dsdx.$$
Novikov's condition and Girsanov's theorem for white noise (see Corollary 3.1.3 in \cite{HAd})
imply that, for $n\in\N$, $\Wtmn=\{\WttBn,\Ft;0\le t\le T,B\in\BR\}$ is a white noise stopped at time $\taunU$, under the
probability measure $\Ptn$ defined on $\FT$ by the recipe
\begin{equation}
\begin{split}
&\frac{d\Ptn}{d\P}=\RNfTtaunW
\\&\eqdef\ \exp\left[\intop_{[0,T\wedge\taunU]\times\R}\RU(s,x)\right.\Wm^{(2)}(ds,dx)
\left.- \frac12\intop_{[0,T\wedge\taunU]\times\R} \RUsq(s,x) ds dx\right].
\label{RNf}
\end{split}
\end{equation}
It then follows that $\solUWt$, $\filTpspace$ is a solution to the $\ehNabdh$ on $\RtaunT$
$\eqdef[0,T\wedge\taunU]\times\R$,
for each $n\in\N$.  Consequently for an arbitrary set $\Lambda\in{\scr B(\CRTR)}$ we get
\begin{equation}
\begin{split}
\Q[V(\cdot,\cdot)\in\Lambda,\taunV=T]&=\Ptn[U(\cdot,\cdot)\in\Lambda,\taunU=T]
\\&=\EP\left[1_{\{U(\cdot,\cdot)\in\Lambda,\taunU=T\}}\RNfTtaunW\right];\quad n\in\N.
\end{split}
\label{1}
\end{equation}
To see \eqref{1} observe that, on the event $\Omega_n^U\eqdef\{\omega\in\Omega^{(2)};\taunU(\omega)=T\}$,
$\solUWt$ is a solution to $\ehNabdh$ on $\RT$, under $\Ptn$, and so the uniqueness in law for $\ehNabdh$ and the
definitions of $\taunU$ and $\taunV$ give the first equality in \eqref{1}.
By the $L^2$ assumption on $\RV$ and the definition of $\taunV$, we have $\lim_{n\to\infty} \Q[\taunV=T]=1$ so that taking
limits in \eqref{1} we get
\begin{equation}
\Q[V(\cdot,\cdot)\in\Lambda]=\lim_{n\to\infty}\Ptn[U(\cdot,\cdot)\in\Lambda,\taunU=T]
=\lim_{n\to\infty}\EP\left[1_{\{U(\cdot,\cdot)\in\Lambda,\taunU=T\}}\RNfTtaunW\right].
\label{2}
\end{equation}
Clearly, if $\P[U(\cdot,\cdot)\in\Lambda]=0$ then $\EP\left[1_{\{U(\cdot,\cdot)\in\Lambda,\taunU=T\}}\RNfTtaunW\right]=0$
for each $n$, and so
$$\Q[V(\cdot,\cdot)\in\Lambda]=\lim_{n\to\infty}\EP\left[1_{\{U(\cdot,\cdot)\in\Lambda,\taunU=T\}}\RNfTtaunW\right]=0.$$
I.e., $\LawVQ$ is absolutely continuous with respect to $\LawUP$ $($on $\scr B(\CRTR)$.    A similar argument yields the
absolute continuity of  $\LawUP$ with respect to $\LawVQ$, and we will omit it.
\end{pf*}
Our compact support result for the Allen-Cahn SPDE \eqnref{AC} can now be proved.
\vspace{2mm} \\
\begin{pf*}{Proof of \thmref{CSPU}}
To see the compact support transfer among \eqref{Pb} and \eqref{Pbd}, observe that if
$\P[U(\cdot,\cdot)\in$ $C_c(\RT;\R)]=1$ ($\P[U(t,\cdot)\in C_c(\R;\R) ]=1$), then by
\thmref{AC1} and \corref{AC2} we have $\Q[V(\cdot,\cdot)\in C_c(\RT;\R) ]=1$ ($\Q[V(t,\cdot)$ $\in C_c(\R;\R)$$]=1$),
respectively, and vice versa.

If $\frac12\le\gamma<1$ and the integrability condition \eqref{integ} is satisfied by solutions of
the Allen-Cahn SPDE \eqref{AC}, then by \corref{ACandheat} the law of \eqref{AC} is equivalent to that
of the associated heat SPDE (without the Allen-Cahn nonlinearity).  Now, observe that
if $U$ is a solution to the heat SPDE associated with \eqref{AC}; then by
\cite{Kr,Mu} we have that, for each $t\in[0,T]$, $U(t,\cdot)$ has compact support (in the space variable) almost surely if
$h(x)$ has compact support and if $0<\gamma<1$.   It then follows, as in the proof of the first part of \thmref{CSPU}
(the compact support transfer among \eqref{Pb} and \eqref{Pbd}), that if $V$ is
a solution to the Allen-Cahn SPDE \eqref{AC}; then, for each $t\in[0,T]$,  $V(t,\cdot)$ has compact support (in space) almost surely
whenever $h(x)$ is compactly supported and $\frac12\le\gamma<1$.  In the opposite direction, the compact
supportedness of $V(t,\cdot)$ for each $t\in[0,T]$ trivially implies the integrability in \eqnref{integ} for $\frac12\le\gamma<1$.
\end{pf*}


\begin{thebibliography}{10}
\bibitem{AllenC} {Allen, S. and Cahn, J.~(1979)}. A microscopic theory for antiphase boundary motion and its
application to antiphase domain coarsing. Acta Metall. {27} 1084--1095.
 \bibitem{HAd} {Allouba, H.~(1996)}.  Different types of SPDEs: existence, uniqueness, and Girsanov theorem.
Ph.D. Dissertation, Cornell University.
 \bibitem{HA1} {Allouba, H.~(1998)}. Different types of SPDEs in the eyes of Girsanov theorem.
 Stochastic Anal.~Appl. {16}, no. 5, 787--810.
 \bibitem{HA2} {Allouba, H.~(1998)}. A non-nonstandard proof of Reimers' existence result for heat SPDEs
J. Appl. Math. Stoch. Anal. {11}, no. 1, 29--41.
\bibitem{Allouba1} {Allouba, H.~(2000)}. Uniqueness in law for the Allen-Cahn SPDE via change of measure. C.R. Acad.~Sci.
 330,  no. 5, 371-376.
\bibitem{HA} {Allouba, H.~(2000)}. SDDEs limits solutions to generalized KPP and Allen-Cahn SPDEs. In preparation.
\bibitem{KM}{El Karoui, N. and M\'el\'eard, S.~(1990)}.  Martingale Measures and Stochastic Calculus.
Probab. Theory Related Fields. {84} 83--101.
\bibitem{DaP}{Da Pratto, G. and Zabczyk, J.~(1992)}. Stochastic equations in infinite dimensions. Cambridge
University Press, Cambridge.
\bibitem{Shreve} {Karatzas, I. and Shreve, S.~(1988)}.  Brownian motion and stochastic calculus. Springer, New York.
\bibitem{Kr} {Krylov, N.~(1997)}. On a result of C. Mueller and E. Perkins.  Probab. Theory Related Fields
{108}, no. 4, 543--557.
\bibitem{Mu1} { Mueller, C.~(1991)}. On the support of solutions to the heat equation with noise.
Stochastics {37}, no. 4, 225--245.
\bibitem{Mu} { Mueller, C. and Perkins, E~(1992)}. The compact support property
for solutions to the heat equation with noise. Probab. Theory Related Fields { 93}, no. 3, 325--358.
\bibitem{Yor} {Revuz, D. and Yor, M.~(1991)}. Continuous martingales
and Brownian motion. Springer, New York.
\bibitem{RC}{Roelly-Coppoletta, S.~(1986)}.  A criterion of convergence of measure-valued processes: application
to measure branching processes. Stochastics { 17}, no. 1-2, 43--65.
\bibitem{WA} {Walsh, J. B.~(1986)}. An introduction to stochastic partial
differential equations. Ecole d'Et\'e de Probabilit\'es de Saint-Flour
XIV.Lecture Notes in Math. {1180}. Springer, New York.
\end{thebibliography}
\end{document}